\documentclass[10pt]{article}
\usepackage{amsmath}
\usepackage{amsfonts}
\usepackage{amssymb} 
\usepackage{authblk}
\usepackage[numbers]{natbib}
\usepackage{geometry}
\usepackage{graphicx}
\usepackage{xcolor}
\usepackage{amssymb}
\usepackage{amsmath}
\usepackage{mathtools}
\usepackage[utf8]{inputenc}
\usepackage{amsthm}
\usepackage{bm}
\usepackage{xcolor}

\usepackage{amsthm}
\usepackage{booktabs} 

\DeclareMathOperator{\mna}{MNA}
\DeclareMathOperator{\nas}{NAS}

\DeclareSymbolFont{eulargesymbols}{U}{zeuex}{m}{n}
\DeclareMathSymbol{\intop}{\mathop}{eulargesymbols}{"52}
\DeclareMathSymbol{\ointop}{\mathop}{eulargesymbols}{"49}

\newcommand{\bN}{\mathbb N}

\newcommand{\bR}{\mathbb{R}}

\newcommand{\bZ}{\mathbb Z}

\newcommand{\e}{\mathrm{e}}

\newcommand{\eps}{\varepsilon}

\newcommand{\la}{\lambda}

\newcommand{\nf}{\infty}

\newcommand{\ph}{\varphi}

\newcommand{\si}{\sigma}
\newcommand{\Si}{\Sigma}

\newcommand{\tht}{\theta}

\renewcommand{\arraystretch}{1.2}

\renewcommand{\ge}{\geqslant}

\renewcommand{\i}{\mathrm{i}}

\renewcommand{\le}{\leqslant}

\def\YYint#1#2#3{{\setbox0=\hbox{$#1{#2#3}{\int}$}
\vcenter{\hbox{$#2#3$}}\kern-.52\wd0}}

\numberwithin{equation}{section}

\newtheorem{lemma}{Lemma}[section]
\newtheorem{theorem}[lemma]{Theorem}

\newtheorem{conjecture}[lemma]{Conjecture}

\newtheorem{example}{Example}[section]
\newtheorem{definition}{Definition}[section]
\newtheorem{remark}[lemma]{Remark}

\author[1]{Manuel Bogoya}
\author[2,3]{Stefano Serra-Capizzano}
\author[4]{Paris Vassalos}
\affil[1]{\small{Departamento de Matem\'aticas, Universidad del Valles, Cali, Colombia}}
\affil[2]{Dipartimento di Scienza e Alta Tecnologia, Insubria University,  Como, Italy}
\affil[3]{Department of Information Technology, Division of Scientific Computing, Uppsala University, Uppsala,  Sweden}
\affil[4]{Department of Informatics, Athens University of Economics and Business, Athens, Greece}
\affil[1]{Email: johan.bogoya@correounivalle.edu.co}
\affil[2,3]{Email: s.serracapizzano@uninsubria.it; stefano.serra@it.uu.se}
\affil[4]{Email: pvassal@aueb.gr}

\title{Fast Toeplitz eigenvalue computations, joining interpolation-extrapolation matrix-less algorithms and simple-loop conjectures: the preconditioned setting}
\date{\today}

\providecommand{\keywords}[1]{\textbf{\textit{Keywords---}} #1}
\begin{document}
\maketitle
\begin{abstract}
Under appropriate technical assumptions, the simple-loop theory allows to deduce various types of asymptotic expansions for the eigenvalues of Toeplitz matrices $T_{n}(f)$ generated by a function $f$, unfortunately, such a theory is not available in the preconditioning setting, that is for matrices of the form $T_{n}^{-1}(g)T_{n}(l)$ with $l,g$ real-valued, $g$ nonnnegative and not identically zero almost everywhere. Independently and under the milder hypothesis that $f=\frac{l}{g}$ is even and monotonic over $[0,\pi]$, matrix-less algorithms have been developed for the fast eigenvalue computation of large preconditioned matrices of the type above, within a linear complexity in the matrix order: behind the high efficiency of such algorithms there are the expansions as in the case $g\equiv 1$, combined with the extrapolation idea, and hence we conjecture that the simple-loop theory has to be extended in such a new setting, as the numerics strongly suggest.

Here we focus our attention on a change of variable, followed by the asymptotic expansion of the new variable, and we consider new matrix-less algorithms ad hoc for the current case.

Numerical experiments show a much higher precision till machine precision and the same linear computation cost, when compared with the matrix-less procedures already proposed in the literature.
\end{abstract}


\maketitle

\section{Introduction}

In this work we design and test fast procedures for the computation of all the eigenvalues of large preconditioned Toeplitz matrices of the form $X_{n}\equiv T_{n}^{-1}(g)T_{n}(l)$ with $g,l$ even, real-valued on $Q\equiv(-\pi,\pi)$, $g>0$ on $(0,\pi)$ such that $f\equiv\frac{l}{g}$ is monotone in the interval $(0,\pi)$. For the formal definition of Toeplitz matrix generated by a Lebesgue integrable function over $Q$ see the first lines of Section~\ref{sc:prel}.

Taking into consideration a clear numerical evidence developed in a systematic series of numerical tests, in \cite{EkGa18} the second author formulated the following conjecture.

\begin{conjecture}\label{conj}
Let $l,g$ be two even functions with $g>0$ on $(0,\pi)$, and suppose that $f\equiv \frac{l}{g}$ is monotone increasing over $(0,\pi)$. Set $X_{n}\equiv T_{n}^{-1}(g)T_{n}(l)$ for all $n$. Then, for every integer $K\geq 0$, every $n$ and every $j=1,\ldots,n$, the following asymptotic expansion holds:
\begin{equation}\label{hoapp}
\la_{j}(X_{n})=f(\tht_{j,n})+\sum_{k=1}^{K}c_{k}(\tht_{j,n})h^{k}+E_{j,n,K},
\end{equation}
where
\begin{itemize}
\item the eigenvalues of $X_{n}$ are arranged in non-decreasing order, $\la_{1}(X_{n})\leq\cdots\leq\la_{n}(X_{n})$;\,\footnote{Note that the eigenvalues of $X_{n}$ are real, because $T_{n}(g)$ is symmetric positive definite and $X_{n}$ is similar to the symmetric matrix $T_{n}^{-\frac{1}{2}}(g)T_{n}(l)T_{n}^{-\frac{1}{2}}(g)$.}
\item $\{c_{k}\}_{k=1}^{\nf}$ is a sequence of functions from $(0,\pi)$ to $\bR$ which depends only on $l$ and $g$;
\item $h\equiv\frac{1}{n+1}$ and $\tht_{j,n}\equiv\frac{j\pi}{n+1}=j\pi h$;
\item $E_{j,n,K}=O(h^{K+1})$ is the remainder (the error), which satisfies the inequality $|E_{j,n,K}|\le c h^{K+1}$ for some constant $c$ depending only on $K,l,g$.
\end{itemize}
\end{conjecture}
\noindent As already mentioned, Conjecture~\ref{conj} was originally formulated and supported through numerical experiments in \cite{EkGa18}. Then the algorithmic proposal was extended and refined in \cite{AhAl18,EkFu18a,EkFu18b,EkGa19}. 
When $g\equiv 1$ and $l$ satisfies further technical additional assumptions, those of the simple-loop method, Conjecture~\ref{conj} was formally proved by Bogoya, B\"ottcher, Grudsky, and Maximenko in a series of papers \cite{BoBo15a,BoBo16,BoGr17,BoBo17}. For a positive function $g$, relation (\ref{hoapp}) was proven, using only purely matrix-theoretic tools, and only for $K=1$ in \cite{AhAl18}.

Here we formulate a second conjecture and we exploit it for designing an even more precise method, when compared with that proposed in \cite{AhAl18,EkGa19}, having the same linear complexity. In fact, the distribution results reported in Theorem~\ref{teo-Sergo} and the localization results in Theorem~\ref{teo-Sext} imply that $\la_{j}(X_{n})=f(s_{j,n})$, $X_{n}=T_{n}^{-1}(g)T_{n}(l)$, $f=\frac{l}{g}$ with $s_{j,n}$ belonging to $(0,\pi)$, and with the sequence 
\[
\Big\{\{s_{j,n}\}_{j=1}^{n}\Big\}_{n}
\] 
distributed as the identity function. 

More precisely, from the combination of Theorem~\ref{teo-Sergo},  Theorem~\ref{teo-Sext}, and of the monotonicity of $f$, for every $j,n$ we easily deduce that
\begin{equation}\label{eq:MainExp zero-level}
\la_{j}(X_{n})\equiv f(s_{j,n}), \quad s_{j,n}=\tht_{j,n}+ o(1),\quad \tht_{j,n}\equiv\frac{j\pi}{n+1}.
\end{equation}
However, we think a richer result holds and more in detail we conjecture that, also in the independent variable $\tht$, for any $K>0$, there exists an asymptotic expansion, regarding exactly the points $s_{j,n}$, given explicitly by the following conjecture.

\begin{conjecture}\label{cj:MainExp}
Let $l,g$ be two even functions with $g>0$ on $(0,\pi)$, and suppose that $f\equiv \frac{l}{g}$ is monotone increasing over $(0,\pi)$. Set $X_{n}\equiv T_{n}^{-1}(g)T_{n}(l)$ for all $n$. Then, for some integer $K\ge0$, every $n$ and every $j=1,\ldots,n$, the following asymptotic expansion holds:
\begin{equation*}
\la_{j}(X_{n})\equiv f(s_{j,n}), \quad s_{j,n}=\tht_{j,n}+\sum_{k=1}^{K}\rho_{k}(\tht_{j,n})h^{k}+E_{j,n,K},
\end{equation*}
where
\begin{itemize}
\item the numbers $s_{j,n}$ are arranged in nondecreasing order;
\item $h\equiv\frac{1}{n+1}$ and $\tht_{j,n}\equiv \pi j h$;
\item the coefficients $\rho_{k}$ are continuous functions from $(0,\pi)$ to $\bR$ which depend only on $f$,
\item $E_{j,n,K}=O(h^{K+1})$ is the remainder (the error), which satisfies the inequality $|E_{j,n,K}|\le c h^{K+1}$ for some constant $c$ depending only on $K,l,g$.
\end{itemize}
\end{conjecture}

This article deals with the adaptation of the interpolation-extrapolation algorithms to the previous change of variable, joined with a trick at the end points introduced in \cite{BoSe21}. 

The numerical results are extremely precise, even compared with the already good performances described in \cite{AhAl18,EkFu18a,EkFu18b,EkGa19,EkGa18}
and of the same order as (or even better than) in \cite{BoEk22} for the non preconditioned setting: in fact, it is not difficult to reach machine precision and the complexity for computing all the eigenvalues is still linear.

  The present work is organized as follows. Preliminary definitions, tools, and results are concisely reported in Section~\ref{sc:prel}. Section~\ref{sc:algo} presents the new adapted algorithm for computing the Toeplitz eigenvalues: as in \cite{EkGa19}, our technique combines the extrapolation procedure proposed in \cite{AhAl18,EkGa18} -- which allows the computation of {\em some} of the eigenvalues of $X_{n}$ -- with an appropriate interpolation process, designed for the simultaneous computation of {\em all} the eigenvalues of $X_{n}$, with the additional end point trick in \cite{BoSe21}. In Section~\ref{sc:num} we present the numerical experiments, while in Section~\ref{conclusions} we draw conclusions and we list few open problems for future research lines, to be investigated in the next future.

\section{Preliminaries and Tools}\label{sc:prel}

For a real or complex valued function $f$ in $L^{1}[-\pi,\pi]$, let $\mathfrak{a}_{j}(f)$ be its $j$th Fourier coefficient, i.e.
\[
\mathfrak{a}_{j}(f)\equiv\frac{1}{2\pi}\int_{-\pi}^{\pi} f(\tht)\e^{-\i j\tht}\,\textrm{d}\tht,\quad j\in\bZ,
\]
and consider the sequence $\{T_{n}(f)\}_{n=1}^{\nf}$ of the $n\times n$ Toeplitz matrices defined by $T_{n}(f)\equiv\big(\mathfrak{a}_{j-k}(f)\big)_{j,k=0}^{n-1}$. The function $f$ is customarily referred to as the generating function of this sequence.

As a second step, we introduce some notations and definitions concerning general sequences of matrices. For any function $F$ defined on the complex field and for any matrix $A_{n}$ of size $d_{n}$, by the symbol $\Si_{\la}(F,A_{n})$, we denote the mean
\begin{equation*}
\Si_{\la}(F,A_{n})\equiv\frac{1}{d_{n}} \sum_{j=1}^{d_{n}}
F(\la_{j}(A_{n})),
\end{equation*}

\begin{definition}
Given a sequence $\{A_{n}\}$ of matrices of size $d_{n}$ with $d_{n}<d_{n+1}$, and given a Lebesgue-measurable function $\psi$ defined over a measurable set $X\subset {\bR}^{\nu}$, $\nu \in \bN^+$, of finite and positive Lebesgue measure $\mu(X)$, we say that $\{A_{n}\}$ is distributed as $(\psi,X)$ in the sense of the eigenvalues if, for any continuous $F$ with bounded support, the following limit relation holds
\begin{equation*}
\lim_{n\rightarrow \nf}\Si_{\la}(F,A_{n})= \frac{1}{\mu(X)}\int_{X} F(\psi)\,\mathrm{d}\mu.
\end{equation*}
In this case, we write in short $\{A_{n}\}\sim_{\la} (\psi,X)$. 
\end{definition}

In Remark~\ref{rem:meaning-distribution} we provide an informal meaning of the notion of eigenvalue distribution. 

\begin{remark}\label{rem:meaning-distribution}
The informal meaning behind the above definition is the following. If $\psi$ is continuous, $n$ is large enough, and
\begin{equation*}
\left\{{\bf x}_{j}^{(d_{n})},\ j=1,\ldots, d_{n}\right\}
\end{equation*}
is an equispaced grid on $X$, then a suitable ordering $\la_{j}(A_{n})$, $j=1,\ldots,d_{n}$, of the eigenvalues of $A_{n}$ is such that the pairs $\big\{\big({\bf x}_{j}^{(d_{n})},\la_{j}(A_{n})\big),\ j=1,\ldots,d_{n}\big\}$ reconstruct approximately the hypersurface
\begin{equation*}
\{({\bf x},\psi({\bf x})),\ {\bf x}\in X\}.
\end{equation*}
 In other words, the spectrum of $A_{n}$ `behaves' like a uniform sampling of $\psi$ over $X$.
For instance, if
$\nu=1$, $d_{n}=n$, and $X=[a,b]$, then the eigenvalues of $A_{n}$ are approximately equal to $\psi\big(a+\frac{j}{n+1}(b-a)\big)$, $j=1,\ldots,n$, for $n$ large enough and up to at most $o(n)$ outliers.
Analogously, if we have $\nu=2$, $d_{n}=n^{2}$, and $X=[a_{1},b_{1}]\times [a_{2},b_{2}]$, then the eigenvalues of the matrix $A_{n}$ are approximately equal to $\psi\big(a_{1}+\frac{j}{n+1}(b_{1}-a_{1}),a_{2}+\frac{k}{n+1}(b_{2}-a_{2})\big)$, $j,k=1,\ldots,n$, for $n$ large enough and up to at most $o(n^{2})$ outliers.
\end{remark}

The asymptotic distribution of eigenvalues and singular values of Toeplitz matrix sequences has been studied deeply and continuously in the last century (for example see \cite{BaGa20b,BaGa20a,BoSi99,GaSe17,GaSe18} and references therein).
For the preconditioned matrix sequences as defined before, a formally similar theory holds, studied by the second author in a series of papers 
\cite{DiBFi93,Se94,Se97b,Se99e,Se99f,Se98a}. 

\begin{theorem}{\rm\cite{Se98a}}
\label{teo-Sergo}
If $g,l$ are integrable over $Q=(-\pi,\pi)$, $g$ is nonnegative almost everywhere (a.e) and not identically zero a.e., and $f=\frac{l}{g}$.  
Then, setting $\{X_{n}\}$ the sequence of preconditioned Toeplitz matrices with $X_{n}=T_{n}^{-1}(g)T_{n}(l)$, we deduce
\begin{equation*}
\{X_{n} \}\sim_{\la} (f,\tilde Q),
\end{equation*}
where $\tilde Q$ is defined as $Q$ minus the set where both $g$ and $l$ vanish simultaneously.
\end{theorem}

\begin{theorem}{\rm \cite{Se97b}}
\label{teo-Sext}
If $g,l$ are integrable over $Q=(-\pi,\pi)$, $g$ is nonnegative almost everywhere (a.e) and not identically zero a.e., and $f=\frac{l}{g}$.  
Then, setting $X_{n}=T_{n}^{-1}(g)T_{n}(l)$, $m$ as the essential infimum of $f$, and $M$ as the essential supremum of $f$, with $m<M$, we deduce
\begin{equation*}
\la_{j}(X_{n})\in (m,M),
\end{equation*}
for all $j=1,\ldots n$, and all $n\ge 1$.
If $m=M$ then the result is trivial since $X_{n}=mI_{n}$ with $I_{n}$ being the identity matrix, so that $\la_{j}(X_{n})\equiv m$, for all $j=1,\ldots n$, and for all $n\ge 1$.
\end{theorem}

We notice that Theorem~\ref{teo-Sergo} reduces to the famous Szeg\H{o}~Theorem \cite{GrSz84} when $g\equiv 1$, in its most general version due to Tyrtyshnikov and Zamarashkin  \cite{TyZa98} (see also the work by Tilli \cite{Ti98a} for the extension to the case of matrix-valued generating functions), while Theorem~\ref{teo-Sext} again for  $g\equiv 1$ reduces to the standard localization results for Toeplitz matrices generated by a Lebesgue integrable function. Finally, it is worth stressing that Remark~\ref{rem:meaning-distribution} of course applies also in this context, but no outliers are present (see Theorem~\ref{teo-Sext} and relations (\ref{eq:MainExp zero-level}) for the precise statements).

\section{Algorithmic proposals}\label{sc:algo}

In the work \cite{BoEk22} we used an asymptotic eigenvalue expansion which was based on the simple-loop theory (see for example \cite{BoBo15a,BoBo16,BoGr17} or the nice review \cite{BoBo17}). However, the preconditioned setting has the additional complication of not having a formal supporting result. Thus, our algorithm has the Conjecture~\ref{cj:MainExp} as its theoretical background.

We considered also the algorithms proposed in \cite{BoSe21,EkFu18b,EkGa19,EkGa18}, and produced again an algorithm suited for parallel implementation and that can be called matrix-less, since it does not require to calculate or even to store the objective matrix entries.

For every $n\in\bN$ let $h\equiv\frac{1}{n+1}$ and $\tht_{j,n}\equiv \pi jh$. One of the key details is that the term $\tht_{j,n}$, remains unchanged for different combinations of $j$ and $n$, for instance
\begin{equation}\label{eq:crule}
\tht_{j,n}=\tht_{cj,c(n+1)-1},
\end{equation}
for any constant $c\in\bN$. Thus, if we select $c=2^{k-1}$ $(k\in\bN)$, we will obtain an increasing sequence of matrix sizes, whose eigenvalues will be used in the precomputing phase.

\begin{figure}[ht]
\centering
\includegraphics[width=0.9\textwidth]{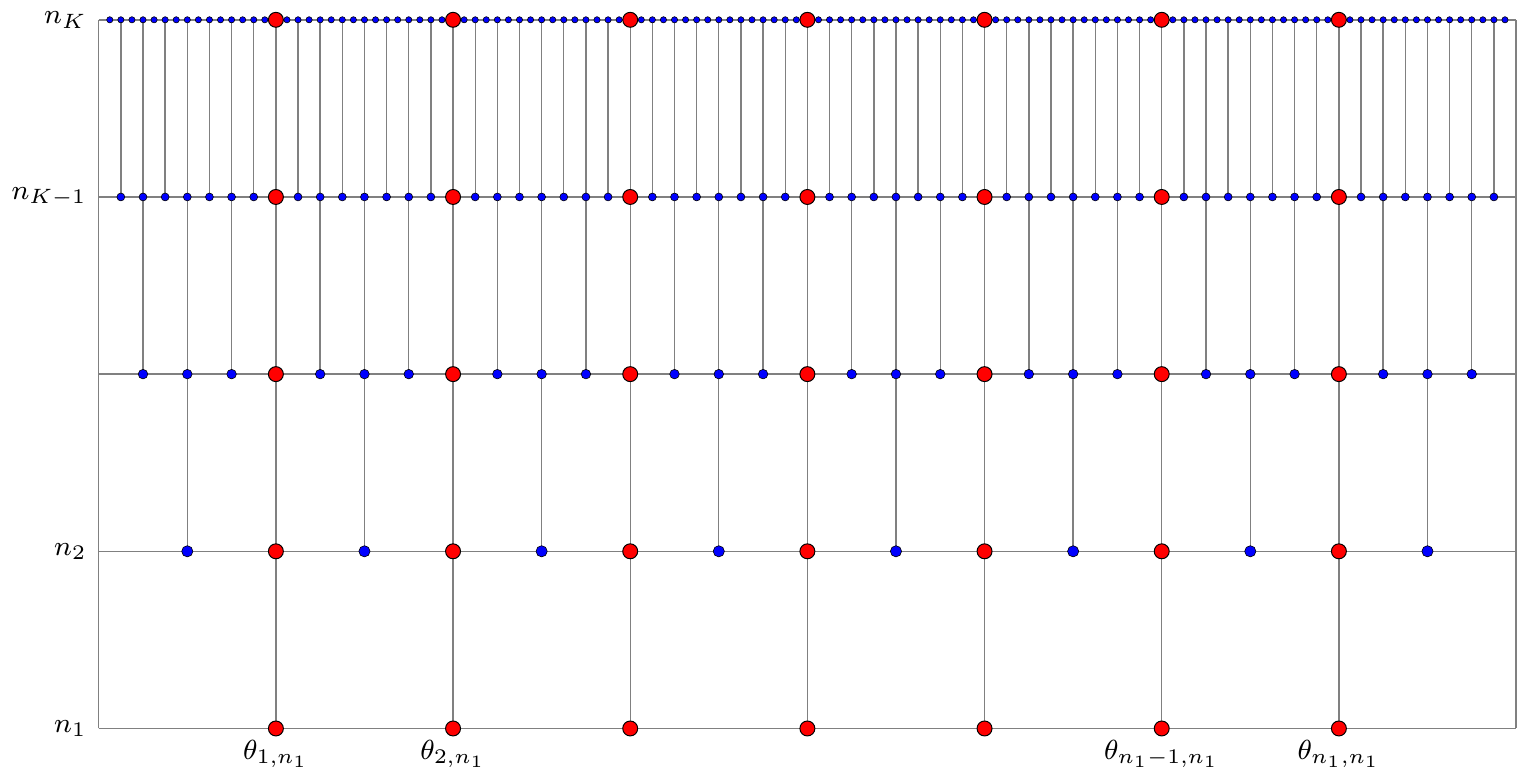}
\caption{The regular grids $\{\tht_{j,n_{k}}\}$ for $j=1,\ldots,n_{k}$ and $k=1,\ldots,K$. In the precomputing phase we will need to calculate the eigenvalues of $X_{n_{k}}$ for $k=1,\ldots,K$ corresponding to the blue and red dots combined, but in the interpolation phase, we will use only the eigenvalues corresponding to the red dots, that is $\la_{j_{k}}(X_{n_{k}})$ for $j_{1}=1,\ldots,n_{1}$ and $k=1,\ldots,K$.}\label{fg:Grid}
\end{figure}

The sequence $\{\tht_{j,n}\}_{j=0}^{n+1}$ is a regular partition of the interval $[0,\pi]$ with step size $\pi h$. We mirror notations for $n$ and $h$, for example, $h_{k}$ means $\frac{1}{n_{k}+1}$, and so on. We assume that
\begin{itemize}
\item the function $f=\frac{l}{g}$ is even and real-valued, strictly increasing in the interval $[0,\pi]$, and $f(0)=0$;
\item $n_{1}$ and $K$ are fixed natural numbers and $n\gg n_{1}$;
\item for $k=1,\ldots,K$ let $n_{k}\equiv 2^{k-1}(n_{1}+1)-1$;
\item for $j_{1}=1,\ldots,n$ and $k=1,\ldots,K$, let $j_{k}\equiv 2^{k-1}j_{1}$.
\end{itemize}
The index $j_{k}$ depends on $j_{1}$, and similarly, the matrix sizes $n_{k}$ depend on $n_{1}$, but for notation simplicity, we suppressed those dependencies.
Following the rule \eqref{eq:crule}, the numbers $j_{k}$ and the matrix sizes $n_{k}$ were calculated in such a way that
\[\si_{j_{1}}\equiv\tht_{j_{1},n_{1}}=\tht_{j_{2},n_{2}}=\cdots=\tht_{j_{K},n_{K}},\]
see Figure~\ref{fg:Grid}.

For well conditioned matrices of sizes in the order $10^{2}$--$10^{4}$ the eigenvalue computation can be done in any modern standard computer. However, in applications like statistical physics and other relevant applied settings we need to handle dimensions of order $10^{8}$--$10^{12}$, a task impossible even for any modern supercomputer. Additionally, the numerical verification of our Conjecture~\ref{cj:MainExp} opens a door for a formal extension of the simple-loop theory to the preconditioned setting. Thus, our aim is to produce an algorithm capable to calculate eigenvalues for `big' matrices, including those extremely large ones. The results show that we can reach machine precision accuracy easily. Recall that $X_{n}\equiv T_{n}^{-1}(g)T_{n}(l)$. As a precomputing phase we need to calculate the eigenvalues of $X_{n}$ for $n=n_{1},\ldots,n_{K}$. This can be easily done with any standard eigensolver (i.e. \texttt{Eigenvalues} in \textsc{Mathematica} or \texttt{eig} in \textsc{Matlab}).

The algorithm has two phases, the first one includes an extrapolation procedure, and the second one is a local interpolation technique.

\noindent{\textbf{Extrapolation.} For each fixed $j_{1}=1,\ldots,n_{1}$ let $\si_{j_{1}}\equiv\tht_{j_{1},n_{1}}=\cdots=\tht_{j_{K},n_{K}}$ (see the red dots in Figure~\ref{fg:Grid}), and apply $K$ times the expansion in the Conjecture~\ref{cj:MainExp} obtaining
\begin{eqnarray*}
s_{j_{1},n_{1}}-\si_{j_{1}}&=&\rho_{1}(\si_{j_{1}})h_{1}+\rho_{2}(\si_{j_{1}})h_{1}^{2}+\cdots+\rho_{K}(\si_{j_{1}})h_{1}^{K}+E_{j_{1},n_{1},K},\\
s_{j_{2},n_{2}}-\si_{j_{1}}&=&\rho_{1}(\si_{j_{1}})h_{2}+\rho_{2}(\si_{j_{1}})h_{2}^{2}+\cdots+\rho_{K}(\si_{j_{1}})h_{2}^{K}+E_{j_{2},n_{2},K},\\
&\vdots&\\
s_{j_{K},n_{K}}-\si_{j_{1}}&=&\rho_{1}(\si_{j_{1}})h_{K}+\rho_{2}(\si_{j_{1}})h_{K}^{2}+\cdots+\rho_{K}(\si_{j_{1}})h_{K}^{K}+E_{j_{K},n_{K},K}.
\end{eqnarray*}
Let $\hat \rho_{k}(\si_{j_{1}})$ be the approximation of $\rho_{k}(\si_{j_{1}})$ obtained by removing all the error terms $E_{j_{k},n_{k},K}$ and solving the resulting linear system:
\begin{equation}\label{eq:IntPhase}
\begin{bmatrix}
h_{1} & h_{1}^{2} & \cdots & h_{1}^{K}\\
h_{2} & h_{2}^{2} & \cdots & h_{2}^{K}\\
\vdots & \vdots & \ddots & \vdots\\
h_{K} & h_{K}^{2} & \cdots & h_{K}^{K}
\end{bmatrix}
\begin{bmatrix}
\hat \rho_{1}(\si_{j_{1}}) \\
\hat \rho_{2}(\si_{j_{1}}) \\
\vdots\\
\hat \rho_{K}(\si_{j_{1}})
\end{bmatrix}
=
\begin{bmatrix}
s_{j_{1},n_{1}} \\
s_{j_{2},n_{2}} \\
\vdots\\
s_{j_{K},n_{K}}
\end{bmatrix}
-\si_{j_{1}}
\begin{bmatrix}
1\\ 1\\ \vdots\\ 1
\end{bmatrix}.
\end{equation}
Let $\phi$ be the inverse function of $f$ restricted to the interval $[0,\pi]$. The value of each $s_{j_{k},n_{k}}$ can be calculated as $\phi(\la_{j_{k}}(X_{n_{k}}))$. If $\phi$ is not available exactly, it can be found numerically with any standard root finder (i.e. \texttt{FindRoot} in \textsc{Mathematica} or \texttt{fzero} in \textsc{Matlab}). As mentioned in previous works, a variant of this extrapolation strategy was first suggested by Albrecht Böttcher in \cite[\S7]{BoBo15a} and is analogous to the Richardson extrapolation employed in the context of Romberg integration \cite[\S3.4]{StBu10}.
\bigskip\\
\noindent\textbf{Interpolation.} For any index $j\in\{1,\ldots,n\}$ and any level $k=1,\ldots,K$, we will estimate $\rho_{k}(\tht_{j,n})$. If $\tht_{j,n}$ coincides with one of the points in the grid $\{\tht_{1,n_{1}},\ldots,\tht_{n_{1},n_{1}}\}$, then we have the approximations $\hat\rho_{k}(\tht_{j,n})$ from the extrapolation phase for free. In any other case, we will do it by interpolating the data
\begin{equation}\label{eq:ExtPhase}
(\tht_{0,n_{k}},\hat \rho_{k}(\tht_{0,n_{k}})),(\tht_{1,n_{k}},\hat \rho_{k}(\tht_{1,n_{k}})),\ldots,(\tht_{n_{k}+1,n_{k}},\hat \rho_{k}(\tht_{n_{k}+1,n_{k}})),
\end{equation}
for $k=1,\ldots,K$, and then evaluating the resulting polynomial at $\tht_{j,n}$. This interpolation can be done in many ways, but to avoid spurious oscillations explained by the Runge phenomenon \cite[p.78]{Da75a}, and following the strategy of the previous works, we decided to do it considering only the $K-k+5$ points in the grid $\{\tht_{0,n_{1}},\ldots,\tht_{n_{1}+1,n_{1}}\}$ which are closest to $\tht_{j,n}$. Those points can be determined uniquely unless $\tht_{j,n}$ is the mid point of two consecutive points in the grid, in which case we can take any of the two possible choices.

Finally, our eigenvalue approximation with $k$ terms, is given by
\begin{equation}\label{eq:NAS}
\la_{j,k}^{\nas}(X_{n})\equiv f\Big(\tht_{j,n}+\sum_{\ell=0}^{k-1} \hat \rho_{\ell}(\tht_{j,n})h^{\ell}\Big),
\end{equation}
where $k=1,\ldots,K$, and NAS stands for ``Numerical Algorithm in the variable $s_{j,k}$''.

\begin{remark}
Since our algorithm is able to reach machine precision accuracy, in the precomputing phase, we advise to calculate $s_{j_{k},n_{k}}$ with a significant number of precision digits, let's say $60$. In all of our examples we used $K=5$ and $n_{1}=100$, and the precomputing phase was carried in a standard computer in only a few minutes.
\end{remark}
\begin{remark}
According to \cite[Th.3.2]{BoBo15b}, for any $u\in(0,1)$, we have
\[\lim_{n\to\nf}\la_{\lceil un\rceil}(X_{n})=f(\pi u),\]
thus, taking $u=\frac{j_{1}}{n_{1}+1}=j_{1}h_{1}$, we can see that for every $j_{1}=1,\ldots,n_{1}$, the sequences $\big(\la_{j_{k}}(X_{n_{k}})\big)_{k\ge1}$ converge to $f(\tht_{j_{1},n_{1}})$. Consequently, the sequences $\big(s_{j_{k},n_{k}}\big)_{k\ge1}$ converge to $\tht_{j_{1},n_{1}}$. The previous result is a direct consequence of the famous Avram--Parter theorem (see \cite{Av88,Pa86} or the beautiful paper \cite{Ty96}) and agrees with our Conjecture~\ref{cj:MainExp}.
\end{remark}

\section{Numerical evidences}\label{sc:num}

In this section we want to compare our algorithm with the one introduced in \cite[\S4]{BoSe21}. To this aim, we start this section analyzing three examples from \cite[\S3]{EkGa18}, involving Real Cosine Trigonometric Polynomials (RCTP). Our numerical experiments was performed with \textsc{Mathematica v.12 (64~bit)} on a platform with 16GB RAM, using an Intel processor QuadCore IntelCore i7 2.6 GHz.

Recall that $X_{n}\equiv T_{n}^{-1}(g)T_{n}(l)$, $\la_{j,k}^{\nas}$ from \eqref{eq:NAS}, and let $\la_{j,k}^{\mna}(X_{n})$ be the $k$th term approximation of $\la_{j}(X_{n})$, given by the Modified Numerical Algorithm \cite[\S4]{BoSe21}. We use the following notation for the absolute individual errors
\[\eps^{\nas}_{j,n,k}\equiv|\la_{j}(X_{n})-\la_{j,k}^{\nas}(X_{n})|,\qquad
\eps^{\mna}_{j,n,k}\equiv|\la_{j}(X_{n})-\la_{j,k}^{\mna}(X_{n})|,\]
and the respective maximum absolute errors
\[\eps^{\nas}_{n,k}\equiv\max\{\eps^{\nas}_{j,n,k}\colon j=1,\ldots,n\},\qquad
\eps^{\mna}_{n,k}\equiv\max\{\eps^{\mna}_{j,n,k}\colon j=1,\ldots,n\}.\]

\begin{example}\label{ex:1}
Consider the two RCTPs
\begin{equation*}
l(\tht)=2-\cos(\tht)-\cos(2\tht),\qquad
g(\tht)=3+2\cos(\tht),
\end{equation*}
and let $f=\frac{l}{g}$. In this case the function $f$ can be simplified to $f(\tht)=1-\cos(\tht)$ which satisfies our hypothesis and gives us an exact inverse function, i.e. $f^{-1}(\ph)=\arccos(1-\ph)$, $\ph\in(0,2)$. The $j$th Fourier coefficient of $\cos(k\tht)$ is $\frac{1}{2}$ for $j=\pm k$, and $0$ in any other case, thus the Fourier coefficients of $l$ and $g$ can be obtained easily by linearity. For a matrix $X_{n}$ of size $n=4096$, the Figure~\ref{fg:Erru} show the log scale of the individual errors $\eps_{j,n,k}^{\nas}$ and $\eps_{j,n,k}^{\mna}$ for different levels $k$, and the Table~\ref{tb:Erru} show the maximum absolute errors $\eps_{n,k}^{\nas}$ and $\eps_{n,k}^{\mna}$ for different matrix sizes $n$ and different levels $k$.

In this case the matrices $T_{n}(l)$ and $T_{n}(g)$ are banded but $X_{n}$ is not Toeplitz, but dense with exponentially decaying entries. Hence, if we work directly on $X_n$, we lose a lot of information and, as a consequence, the exact computation of its eigenvalues is quite hard. For example, when calculating the eigenvalues of $X_{512}$ with machine precision accuracy, we only get $2$ correct digits. Nevertheless, even for this small matrix, our algorithm is able to produce 15 correct digits and faster.

{\renewcommand{\arraystretch}{1.2}
\begin{table}[ht]
\caption{Example~\ref{ex:1}: The maximum errors $\eps_{n,k}^{\mna}$, $\eps_{n,k}^{\nas}$, and maximum normalized errors $(n+1)^{k}\,\eps_{n,k}^{\nas}$ for the levels $k=1,\ldots,5$, and different matrix sizes $n$, corresponding to the matrix $X_{n}=T_{n}^{-1}(g)T_{n}(l)$ where $l(\tht)=2-\cos(\tht)-\cos(2\tht)$ and
$g(\tht)=3+2\cos(\tht)$. We used a grid of size $n_{1}=100$.}\label{tb:Erru}
\centering
\begin{tabular}{rlllll}
\toprule
\multicolumn{1}{c}{$n$} & \multicolumn{1}{c}{$256$} & \multicolumn{1}{c}{$512$} & \multicolumn{1}{c}{$1024$} & \multicolumn{1}{c}{$2048$} & \multicolumn{1}{c}{$4096$} \\ \midrule

$\eps_{n,1}^{\mna}$ & $2.935\times10^{-3}$ & $1.4706\times10^{-3}$ & $7.3605\times10^{-4}$ & $3.6822\times10^{-4}$ & $1.8416\times10^{-4}$ \\ \midrule
$\eps_{n,1}^{\nas}$ & $2.935\times10^{-3}$ & $1.4706\times10^{-3}$ & $7.3605\times10^{-4}$ & $3.6822\times10^{-4}$ & $1.8416\times10^{-4}$ \\ \midrule
$(n+1)\,\eps_{n,1}^{\nas}$ & $7.5429\times10^{-1}$ & $7.5440\times10^{-1}$ & $7.5445\times10^{-1}$ & $7.5448\times10^{-1}$ & $7.5450\times10^{-1}$ \\ \midrule

$\eps_{n,2}^{\mna}$ & $6.0234\times10^{-6}$ & $1.5189\times10^{-6}$ & $3.8421\times10^{-7}$ & $9.8046\times10^{-8}$ & $2.5538\times10^{-8}$ \\ \midrule
$\eps_{n,2}^{\nas}$ & $3.4682\times10^{-6}$ & $8.6926\times10^{-7}$ & $2.1759\times10^{-7}$ & $5.4432\times10^{-8}$ & $1.3612\times10^{-8}$ \\ \midrule
$(n+1)^{2}\,\eps_{n,2}^{\nas}$ & $2.2907\times10^{-1}$ & $2.2876\times10^{-1}$ & $2.2861\times10^{-1}$ & $2.2853\times10^{-1}$ & $2.2849\times10^{-1}$ \\ \midrule

$\eps_{n,3}^{\mna}$ & $1.8060\times10^{-8}$ & $2.2864\times10^{-9}$ & $8.4778\times10^{-10}$ & $4.2540\times10^{-10}$ & $2.1313\times10^{-10}$ \\ \midrule
$\eps_{n,3}^{\nas}$ & $1.4429\times10^{-8}$ & $1.8129\times10^{-9}$ & $2.2720\times10^{-10}$ & $2.8437\times10^{-11}$ & $3.5569\times10^{-12}$ \\ \midrule
$(n+1)^{3}\,\eps_{n,3}^{\nas}$ & $2.4492\times10^{-1}$ & $2.4476\times10^{-1}$ & $2.4467\times10^{-1}$ & $2.4463\times10^{-1}$ & $2.4461\times10^{-1}$ \\ \midrule

$\eps_{n,4}^{\mna}$ & $1.5184\times10^{-10}$ & $5.7689\times10^{-11}$ & $2.4437\times10^{-11}$ & $1.2058\times10^{-11}$ & $6.0583\times10^{-12}$ \\ \midrule
$\eps_{n,4}^{\nas}$ & $4.9519\times10^{-11}$ & $3.1141\times10^{-12}$ & $1.9522\times10^{-13}$ & $1.2221\times10^{-14}$ & $7.6657\times10^{-16}$ \\ \midrule
$(n+1)^{4}\,\eps_{n,4}^{\nas}$ & $2.1603\times10^{-1}$ & $2.1568\times10^{-1}$ & $2.1548\times10^{-1}$ & $2.1541\times10^{-1}$ & $2.1598\times10^{-1}$ \\ \midrule

$\eps_{n,5}^{\mna}$ & $2.8044\times10^{-11}$ & $1.5290\times10^{-11}$ & $7.9990\times10^{-12}$ & $4.0993\times10^{-12}$ & $2.0649\times10^{-12}$ \\ \midrule
$\eps_{n,5}^{\nas}$ & $1.8256\times10^{-13}$ & $5.7554\times10^{-15}$ & $1.8077\times10^{-16}$ & $5.6588\times10^{-18}$ & $2.3660\times10^{-18}$ \\ \midrule
$(n+1)^{5}\,\eps_{n,5}^{\nas}$ & $2.0467\times10^{-1}$ & $2.0448\times10^{-1}$ & $2.0453\times10^{-1}$ & $2.0438\times10^{-1}$ & $2.7311\times10^{\,0}$ \\ \bottomrule
\end{tabular}
\end{table}}

\begin{figure}[ht]
\centering
\includegraphics[width=0.8\textwidth]{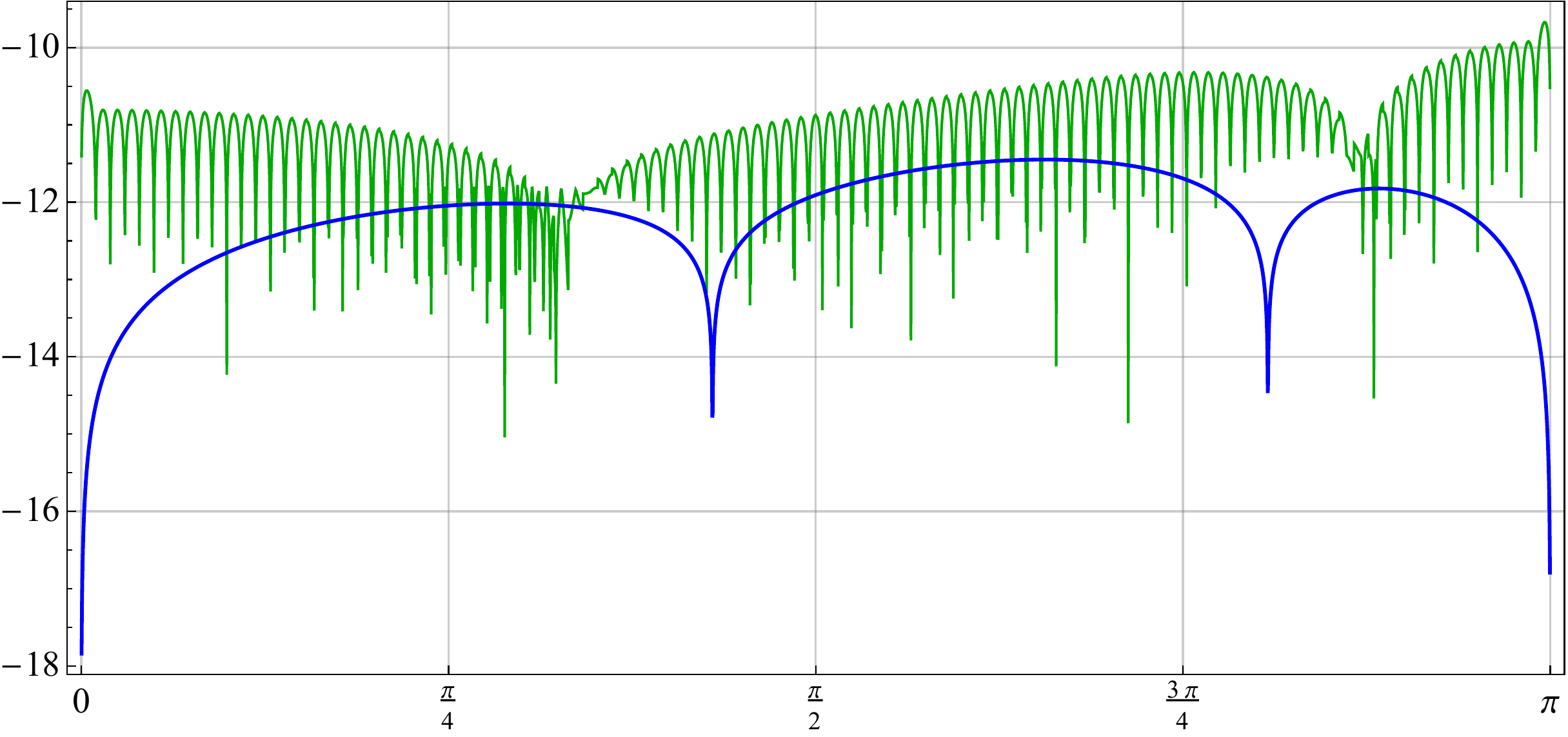}
\put(-189,180){\footnotesize $k=3$}\medskip\\
\includegraphics[width=0.8\textwidth]{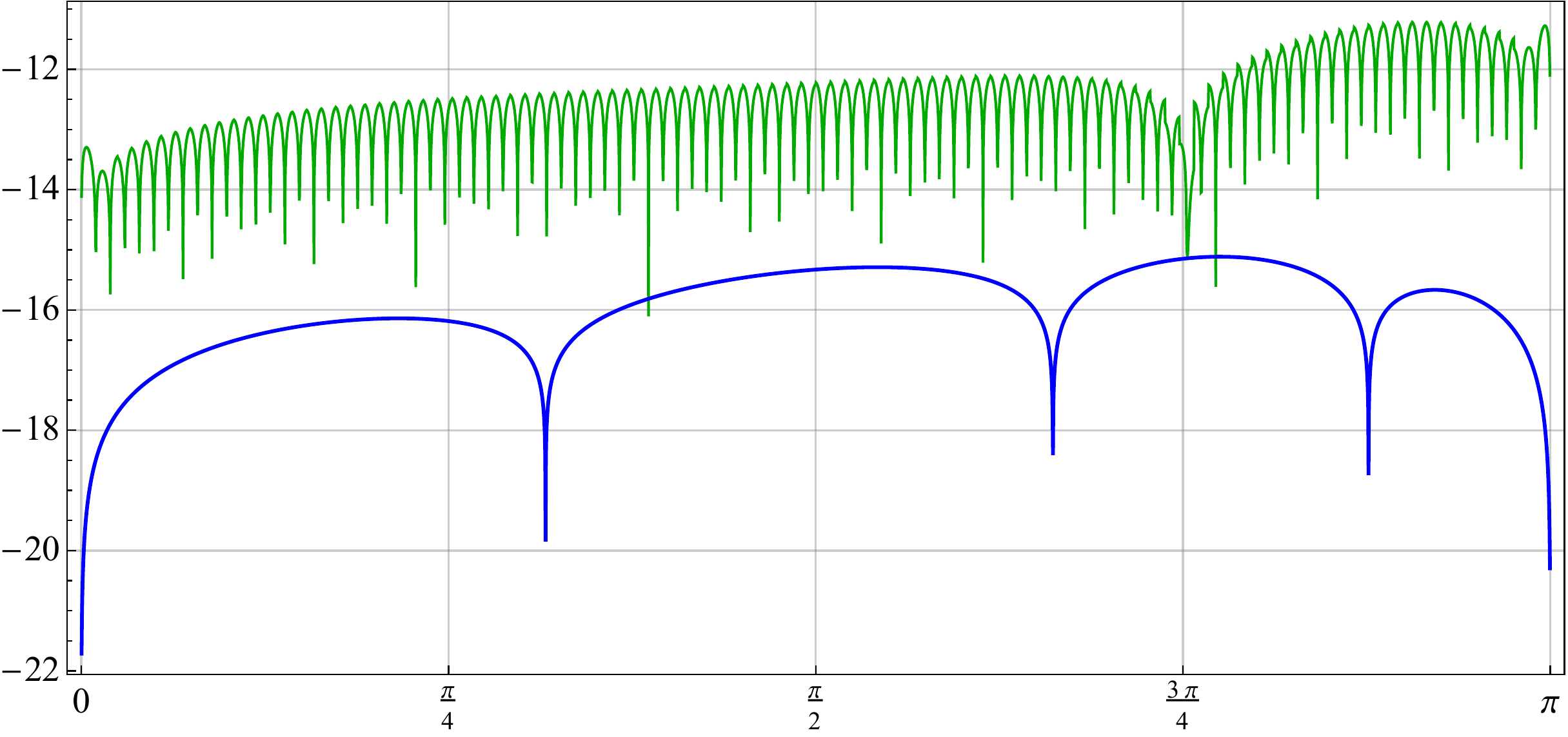}
\put(-189,180){\footnotesize $k=4$}\medskip\\
\includegraphics[width=0.8\textwidth]{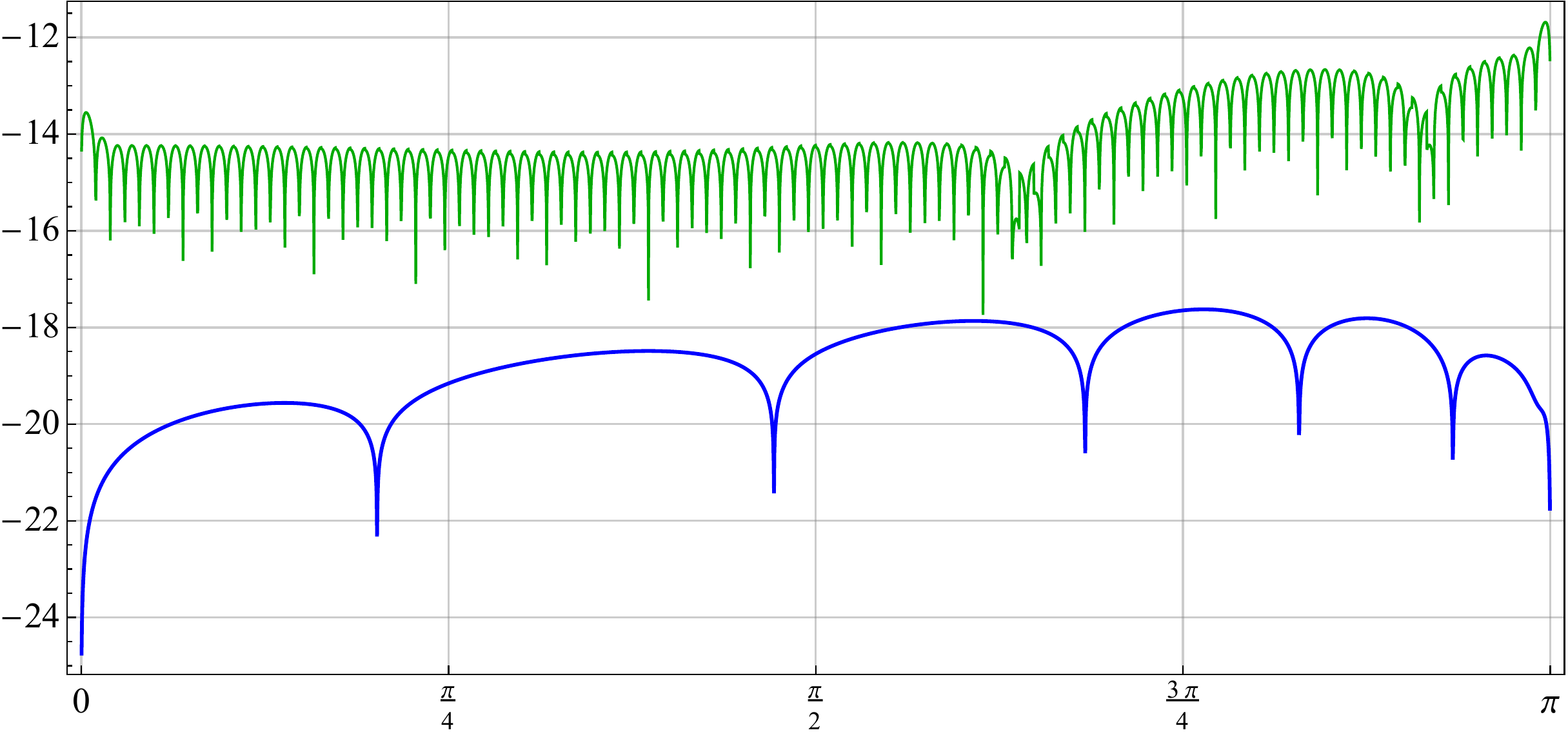}
\put(-189,180){\footnotesize $k=5$}
\caption{Example~\ref{ex:1}: The base-10 logarithm for the individual errors $\eps_{j,n,k}^{\nas}$ (blue) and $\eps_{j,n,k}^{\mna}$ (green) for the matrix $X_{n}=T_{n}^{-1}(g)T_{n}(l)$ with $l(\tht)=2-\cos(\tht)-\cos(2\tht)$ and
$g(\tht)=3+2\cos(\tht)$, a matrix size $n=4096$, a grid size $n_{1}=100$, and different levels $k$.}\label{fg:Erru}
\end{figure}
\clearpage
In \cite[Th.3]{EkGa18} the authors provide an error estimate of the kind
\[\eps_{j,n,K}^{\rm NA}\le c h_{1}^{K}h,\]
where $c$ is a constant depending only on $K,l,g$. The Table~\ref{tb:Erru} shows that the boundary modification introduced in \cite[\S4]{BoSe21}, which we can see reflected in the maximum errors $\eps_{n,k}^{\mna}$, is producing much better results, more specifically, for a grid size of $n_{1}=100$ and similar matrix sizes, they found a maximum error of $\approx 10^{-10}$ while the $\mna$ algorithm produced $\approx10^{-12}$, explaining why they were observing the constant $c$ to grow quickly with $K$. We can also see that our error estimate has the bounding $\eps_{n,k}^{\nas}\le ch^{k}$, for some constant $c$ depending only on $K,l,g$, in perfect agreement with the Conjecture~\ref{cj:MainExp}, and proving that our algorithm is matching the theoretical error and reaching its maximum possible accuracy.
\end{example}

\begin{example}\label{ex:2}
Consider the two RCTPs
\begin{eqnarray*}
l(\tht)&=&40-15\cos(\tht)-24\cos(2\tht)-\cos(3\tht),\\
g(\tht)&=&1208+1191\cos(\tht)+120\cos(2\tht)+\cos(3\tht),
\end{eqnarray*}
and let $f=\frac{l}{g}$. In this case the function $f$ has no important simplification and we must use a numerical root-solver to get its inverse function. As in the Example~\ref{ex:1}, the Fourier coefficients of the symbols $l,g$ can be exactly calculated by linearity and the rule mentioned there. For a matrix $X_{n}$ of size $n=4096$, the Figure~\ref{fg:Errv} show the log scale of the individual errors $\eps_{j,n,k}^{\nas}$ and $\eps_{j,n,k}^{\mna}$ for different levels $k$, and the Table~\ref{tb:Errv} show the maximum absolute errors $\eps_{n,k}^{\nas}$ and $\eps_{n,k}^{\mna}$ for different matrix sizes $n$ and different levels $k$.

{\renewcommand{\arraystretch}{1.2}
\begin{table}
\centering
\caption{Example~\ref{ex:2}: The maximum errors $\eps_{n,k}^{\mna}$, $\eps_{n,k}^{\nas}$, and maximum normalized errors $(n+1)^{k}\,\eps_{n,k}^{\nas}$ for the levels $k=1,\ldots,5$ and different matrix sizes $n$, corresponding to the matrix $X_{n}=T_{n}^{-1}(g)T_{n}(l)$ where $l(\tht)=40-15\cos(\tht)-24\cos(2\tht)-\cos(3\tht)$ and $g(\tht)=1208+1191\cos(\tht)+120\cos(2\tht)+\cos(3\tht)$. We used a grid of size $n_{1}=100$.}\label{tb:Errv}
\begin{tabular}{rlllll}
\toprule
\multicolumn{1}{c}{$n$} & \multicolumn{1}{c}{$256$} & \multicolumn{1}{c}{$512$} & \multicolumn{1}{c}{$1024$} & \multicolumn{1}{c}{$2048$} & \multicolumn{1}{c}{$4096$} \\ \midrule

$\eps_{n,1}^{\mna}$ & $2.9350\times10^{-3}$ & $1.4706\times10^{-3}$ & $7.3605\times10^{-4}$ & $3.6822\times10^{-4}$ & $1.8416\times10^{-4}$ \\ \midrule
$\eps_{n,1}^{\nas}$ & $2.9350\times10^{-3}$ & $1.4706\times10^{-3}$ & $7.3605\times10^{-4}$ & $3.6822\times10^{-4}$ & $1.8416\times10^{-4}$ \\ \midrule
$(n+1)\,\eps_{n,1}^{\nas}$ & $7.5429\times10^{-1}$ & $7.5440\times10^{-1}$ & $7.5445\times10^{-1}$ & $7.5448\times10^{-1}$ & $7.5450\times10^{-1}$ \\ \midrule

$\eps_{n,2}^{\mna}$ & $6.0234\times10^{-6}$ & $1.5189\times10^{-6}$ & $3.8421\times10^{-7}$ & $9.8046\times10^{-8}$ & $2.5538\times10^{-8}$ \\ \midrule
$\eps_{n,2}^{\nas}$ & $3.4682\times10^{-6}$ & $8.6926\times10^{-7}$ & $2.1759\times10^{-7}$ & $5.4432\times10^{-8}$ & $1.3612\times10^{-8}$ \\ \midrule
$(n+1)^{2}\,\eps_{n,2}^{\nas}$ & $2.2907\times10^{-1}$ & $2.2876\times10^{-1}$ & $2.2861\times10^{-1}$ & $2.2853\times10^{-1}$ & $2.2849\times10^{-1}$ \\ \midrule

$\eps_{n,3}^{\mna}$ & $1.8060\times10^{-8}$ & $2.2864\times10^{-9}$ & $8.4778\times10^{-10}$ & $4.2540\times10^{-10}$ & $2.1313\times10^{-10}$ \\ \midrule
$\eps_{n,3}^{\nas}$ & $1.4429\times10^{-8}$ & $1.8129\times10^{-9}$ & $2.2720\times10^{-10}$ & $2.8437\times10^{-11}$ & $3.5569\times10^{-12}$ \\ \midrule
$(n+1)^{3}\,\eps_{n,3}^{\nas}$ & $2.4492\times10^{-1}$ & $2.4476\times10^{-1}$ & $2.4467\times10^{-1}$ & $2.4463\times10^{-1}$ & $2.4461\times10^{-1}$ \\ \midrule

$\eps_{n,4}^{\mna}$ & $1.5184\times10^{-10}$ & $5.7689\times10^{-11}$ & $2.4437\times10^{-11}$ & $1.2058\times10^{-11}$ & $6.0583\times10^{-12}$ \\ \midrule
$\eps_{n,4}^{\nas}$ & $4.9519\times10^{-11}$ & $3.1141\times10^{-12}$ & $1.9522\times10^{-13}$ & $1.2221\times10^{-14}$ & $7.6657\times10^{-16}$ \\ \midrule
$(n+1)^{4}\,\eps_{n,4}^{\nas}$ & $2.1603\times10^{-1}$ & $2.1568\times10^{-1}$ & $2.1548\times10^{-1}$ & $2.1541\times10^{-1}$ & $2.1598\times10^{-1}$ \\ \midrule

$\eps_{n,5}^{\mna}$ & $2.8044\times10^{-11}$ & $1.5290\times10^{-11}$ & $7.9990\times10^{-12}$ & $4.0983\times10^{-12}$ & $2.0649\times10^{-12}$ \\ \midrule
$\eps_{n,5}^{\nas}$ & $1.8256\times10^{-13}$ & $5.7554\times10^{-15}$ & $1.8077\times10^{-16}$ & $1.6588\times10^{-18}$ & $2.3660\times10^{-18}$ \\ \midrule
$(n+1)^{5}\,\eps_{n,5}^{\nas}$ & $2.0467\times10^{-1}$ & $2.0448\times10^{-1}$ & $2.0453\times10^{-1}$ & $2.0438\times10^{-1}$ & $2.7311\times10^{\,0}$ \\ \bottomrule
\end{tabular}
\end{table}}
\clearpage

\begin{figure}[ht]
\centering
\includegraphics[width=0.8\textwidth]{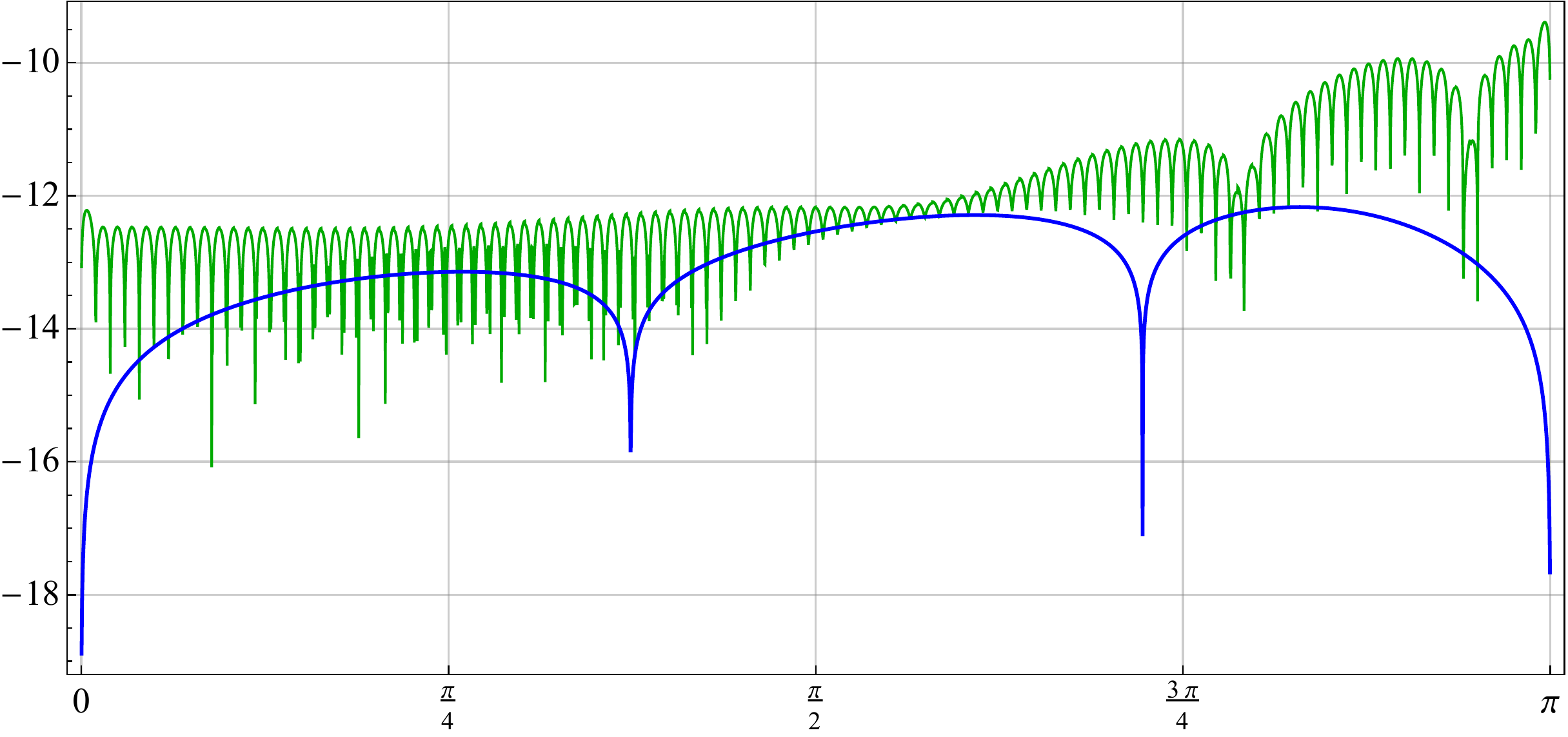}
\put(-189,180){\footnotesize $k=3$}\medskip\\
\includegraphics[width=0.8\textwidth]{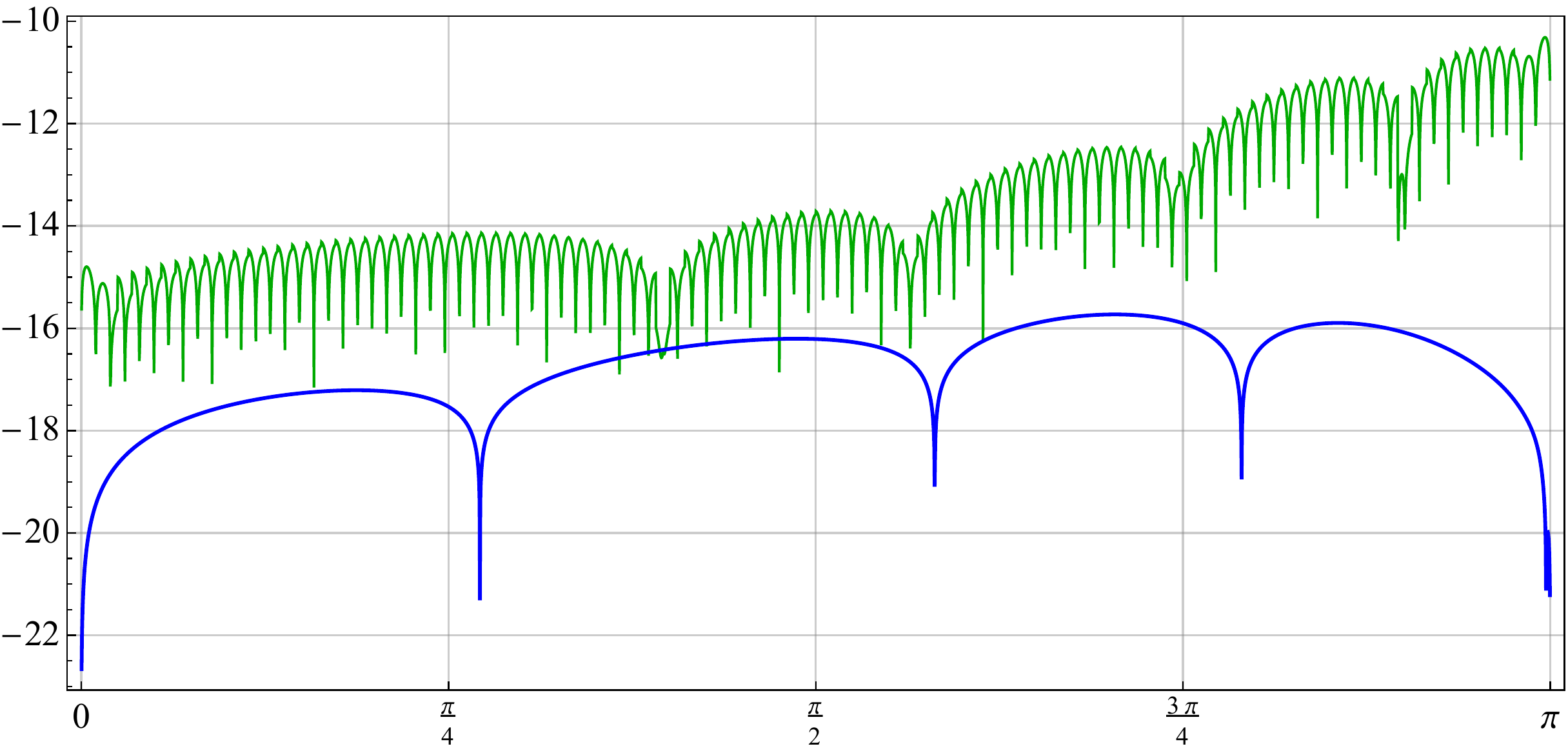}
\put(-189,180){\footnotesize $k=4$}\medskip\\
\includegraphics[width=0.8\textwidth]{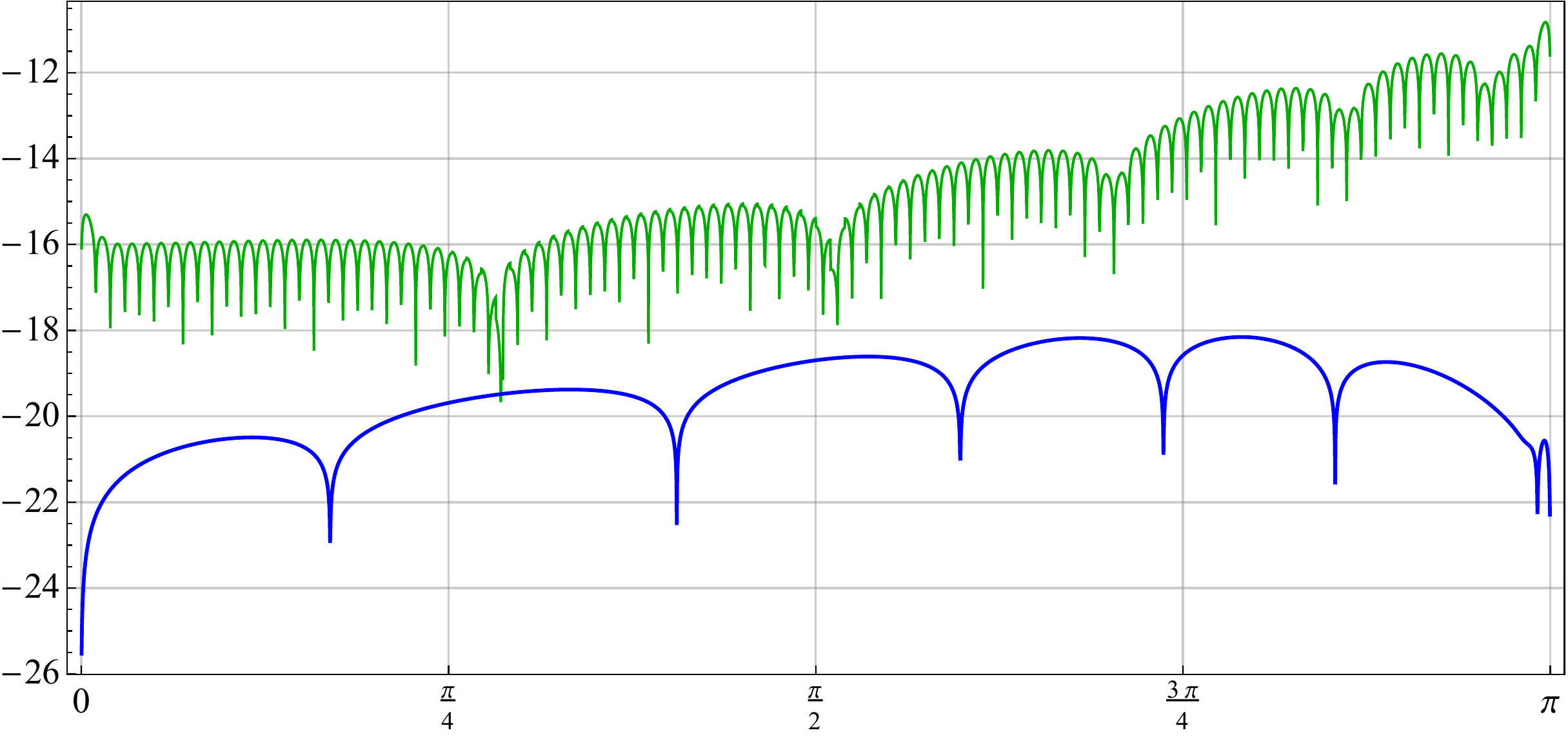}
\put(-189,180){\footnotesize $k=5$}
\caption{Example~\ref{ex:2}: The base-10 logarithm for the individual errors $\eps_{j,n,k}^{\nas}$ (blue) and $\eps_{j,n,k}^{\mna}$ (green) for the matrix $X_{n}=T_{n}^{-1}(g)T_{n}(l)$ with $l(\tht)=40-15\cos(\tht)-24\cos(2\tht)-\cos(3\tht)$ and $g(\tht)=1208+1191\cos(\tht)+120\cos(2\tht)+\cos(3\tht)$, a matrix size $n=4096$, a grid size $n_{1}=100$, and different levels $k$.}\label{fg:Errv}
\end{figure}
\clearpage

There are different alternatives to improve the result of our model, for example, in \cite{EkGa18} the authors took a fixed level $k$ and manage different grid sizes $n_{1}=25,50,100,200,400$, and presented the respective log scaled error figures. Other alternative is to increase the number of levels $K$ in the extrapolation phase \eqref{eq:IntPhase}, or to increase the number of interpolated points in the interpolation phase \eqref{eq:ExtPhase}. In this article, we decided to handle a fixed grid size $n_{1}=100$, $K=5$, $K-k+5$ interpolated points at level $k$, and show the error evolution for different matrix sizes $n$ and different levels $k$. In this way we can observe the algorithm accuracy. As in Example~\ref{ex:1}, the Table~\ref{tb:Errv}, show that our model is nicely following the bound $\eps_{n,k}^{\nas}=O(h^{k})$ which is the maximum possible accuracy.
\end{example}

\begin{example}\label{ex:3}
Consider the two RCTPs
\begin{eqnarray*}
l(\tht)&=&\frac{35}{2}-12\cos(\tht)-6\cos(2\tht)+\frac{1}{2}\cos(4\tht),\\
g(\tht)&=&8-3\cos(\tht)-4\cos(2\tht)-\cos(3\tht),
\end{eqnarray*}
and let $f=\frac{l}{g}$. In this case the function $f$ can be simplified to $f(\tht)=2-\cos(\tht)$ which satisfies our hypothesis and gives us an exact inverse function, i.e. $f^{-1}(\ph)=\arccos(2-\ph)$, $\ph\in(1,3)$. As in the Example~\ref{ex:1}, the Fourier coefficients of the symbols $l,g$ can be exactly calculated by linearity and the rule mentioned there. The Figure~\ref{fg:Errw} and the Table~\ref{tb:Errw} show the data.

\begin{figure}[ht]
\centering
\includegraphics[width=0.8\textwidth]{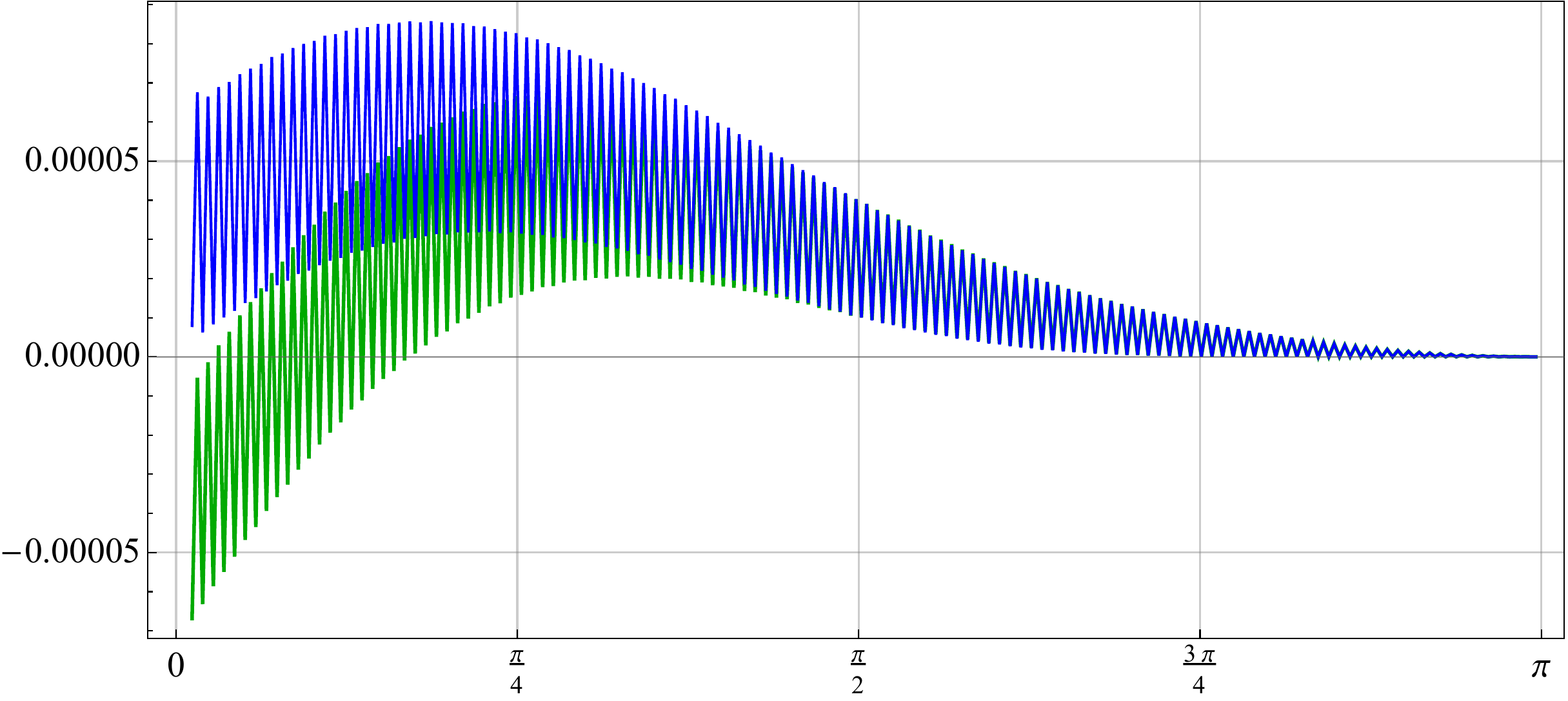}
\caption{Example~\ref{ex:3}: The individual eigenvalue errors $\eps_{j,n,k}^{\nas}$ (blue) and $\eps_{j,n,k}^{\mna}$ (green), for the level $k=2$ and a matrix $X_{n}=T_{n}^{-1}(g)T_{n}(l)$ of size $n=256$. The symbols $l,g$ are given by $l(\tht)=\frac{35}{2}-12\cos(\tht)-6\cos(2\tht)+\frac{1}{2}\cos(4\tht)$ and 
$g(\tht)=8-3\cos(\tht)-4\cos(2\tht)-\cos(3\tht)$. We used a grid of size $n_{1}=100$.}\label{fg:Errw}
\end{figure}

{\renewcommand{\arraystretch}{1.2}
\begin{table}[ht]
\centering
\begin{tabular}{rlllll}
\toprule
\multicolumn{1}{c}{$n$} & \multicolumn{1}{c}{$256$} & \multicolumn{1}{c}{$512$} & \multicolumn{1}{c}{$1024$} & \multicolumn{1}{c}{$2048$} & \multicolumn{1}{c}{$4096$} \\ \midrule

$\eps_{n,1}^{\mna}$ & $4.6910\times10^{-3}$ & $2.3666\times10^{-3}$ & $1.1887\times10^{-3}$ & $5.9570\times10^{-4}$ & $2.9819\times10^{-4}$ \\ \midrule
$\eps_{n,1}^{\nas}$ & $4.6910\times10^{-3}$ & $2.3666\times10^{-3}$ & $1.1887\times10^{-3}$ & $5.9570\times10^{-4}$ & $2.9819\times10^{-4}$ \\ \midrule
$(n+1)\,\eps_{n,1}^{\nas}$ & $1.2056\times10^{\,0}$ & $1.2141\times10^{\,0}$ & $1.2184\times10^{\,0}$ & $1.2206\times10^{\,0}$ & $1.2217\times10^{\,0}$ \\ \midrule

$\eps_{n,2}^{\mna}$ & $6.7245\times10^{-5}$ & $1.7140\times10^{-5}$ & $4.3543\times10^{-6}$ & $1.1152\times10^{-6}$ & $2.9156\times10^{-7}$ \\ \midrule
$\eps_{n,2}^{\nas}$ & $8.5834\times10^{-5}$ & $2.1860\times10^{-5}$ & $6.6691\times10^{-6}$ & $1.0156\times10^{-5}$ & $5.3725\times10^{-6}$ \\ \midrule
$(n+1)^{2}\,\eps_{n,2}^{\nas}$ & $5.6692\times10^{\,0}$ & $5.7530\times10^{\,0}$ & $7.0067\times10^{\,0}$ & $4.2640\times10^{\,1}$ & $9.0180\times10^{1}$ \\ \bottomrule
\end{tabular}
\vspace{2mm}
\caption{Example~\ref{ex:3}: The maximum errors $\eps_{n,k}^{\mna}$, $\eps_{n,k}^{\nas}$, and maximum normalized errors $(n+1)^{k}\,\eps_{n,k}^{\nas}$ for the levels $k=1,2,3$, and different matrix sizes $n$, corresponding to the matrix $X_{n}=T_{n}^{-1}(g)T_{n}(l)$ where the symbols $l,g$ are given by $l(\tht)=\frac{35}{2}-12\cos(\tht)-6\cos(2\tht)+\frac{1}{2}\cos(4\tht)$ and 
$g(\tht)=8-3\cos(\tht)-4\cos(2\tht)-\cos(3\tht)$. We used a grid of size $n_{1}=100$.}\label{tb:Errw}
\end{table}}

The Figure~\ref{fg:Errw} shows that the eigenvalues of the matrix $X_{n}$, with even and odd numbers have different behaviors. This phenomenon was studied in \cite{BaGr17} where the authors deduced formally, asymptotic individual expansions for the eigenvalues of certain penta-diagonal Toeplitz matrix. They produced two different expansions, one for the eigenvalues with even numbers, and another for the odd ones. 

As we can see in the Table~\ref{tb:Errw}, our algorithm was able to produce acceptable results until the level $k=2$ only. We think that it is possible to adapt our algorithm to this case, but it is the topic of a future investigation.
\end{example}

\section{Conclusions}\label{conclusions}

Under appropriate technical assumptions, the simple-loop theory allows to deduce various types of asymptotic expansions for the eigenvalues of Toeplitz matrices $T_{n}(f)$ generated by a function $f$. Independently and under the milder hypothesis that $f$ is even and monotonic over $[0,\pi]$, matrix-less algorithms have been developed for the fast eigenvalue computation of large Toeplitz matrices. These procedures work with a linear complexity in the matrix order $n$ and behind the high efficiency of such algorithms there are the expansions predicted by the simple-loop theory, combined with the extrapolation idea.

Here we conjectured that the same type of expansions hold also for preconditioned matrix sequences, by focusing our attention on a change of variable 
\[
\la_{j}(X_{n})\equiv f(s_{j,n}), \qquad s_{j,n}=\tht_{j,n}+\sum_{k=1}^{K}\rho_{k}(\tht_{j,n})h^{k}+E_{j,n,K},
\]
and then we adapted the matrix-less procedures to the considered new setting.

Numerical experiments have shown, in a clear way, a much higher precision (till machine precision) and the same linear computation cost, when compared with the matrix-less procedures already presented in the relevant literature. 

As next steps the following questions remain to be investigated:

\begin{itemize}
\item Taking inspiration from the works on the simple-loop setting \cite{BoBo15a,BoBo16,BoGr17,BoBo17}, extension of the proofs from the pure Toeplitz setting to the preconditioned Toeplitz setting, as described in this work and as strongly confirmed by the numerical experiments; 
\item Applications to the block cases (see \cite{BaGa20b,BaGa20a} for the theory in the block case) and related applications \cite{EkFu18a,EkFu18b} to differential problems, especially systems of ordinary differential equations (ODEs) \cite{GaMa18}, and/or ODE approximation via Discontinuous Galerkin methods \cite{DuFa18}, Finite Element methods of high order \cite{GaSe15}, Isogeometric Analysis with intermediate smoothness \cite{GaSp19}, etc.;
\item A fine error analysis for having a theoretical explanation of the reason why the new expansion leads in practice to much smaller errors, when compared with the numerical results in \cite{AhAl18,EkGa19}.  
\end{itemize}

\bibliographystyle{plain}
\bibliography{Toeplitz.bib}
\end{document}